\documentclass{amsart}

\usepackage{amssymb}

\newtheorem{thm}{Theorem}%[section]

\newtheorem{lem}[thm]{Lemma}
\newtheorem{cor}[thm]{Corollary}

\newtheorem{prop}[thm]{Proposition}

\theoremstyle{definition}
\newtheorem{defn}[thm]{Definition}

\newtheorem{say}[thm]{}
\newtheorem{exmp}[thm]{Example}

\newtheorem{ack}{Acknowledgments}

\newtheorem{defn-thm}[thm]{Definition--Theorem}  %!!!!!!!!!!!!!!!!!!!!!!!!
\newtheorem{defn-lem}[thm]{Definition--Lemma}  %!!!!!!!!!!!!!!!!!!!!!!!!
\theoremstyle{remark}

\renewcommand{\o}[0]{{\mathcal O}} 
\newcommand{\z}[0]{{\mathbb Z}}

\renewcommand{\a}[0]{{\mathbb A}}

\newcommand{\p}[0]{{\mathbb P}}

\newcommand{\qtq}[1]{\quad\mbox{#1}\quad}
\newcommand{\spec}[0]{\operatorname{Spec}}

\newcommand{\supp}[0]{\operatorname{Supp}}    
\newcommand{\red}[0]{\operatorname{red}}    
\newcommand{\codim}[0]{\operatorname{codim}}    
\newcommand{\im}[0]{\operatorname{im}}

\newcommand{\coker}[0]{\operatorname{coker}}    
    
\newcommand{\Hom}[0]{\operatorname{Hom}}

\newcommand{\aut}[0]{\operatorname{Aut}}

\newcommand{\tors}[0]{\operatorname{tors}}  
\newcommand{\tor}[0]{\operatorname{Tor}}

\newcommand{\hilb}[0]{\operatorname{Hilb}}

\newcommand{\onto}[0]{\twoheadrightarrow}

\newcommand{\ann}[0]{\operatorname{Ann}}

\newcommand{\uhom}[0]{\operatorname{\underline{Hom}}}  
\newcommand{\quot}[0]{\operatorname{Quot}}
  
\newcommand{\mor}[0]{\operatorname{Mor}} 
\newcommand{\depth}[0]{\operatorname{depth}}

\def\into{\DOTSB\lhook\joinrel\to}

\begin{document}
\bibliographystyle{amsalpha}

\title{Hulls and Husks}
\author{J\'anos Koll\'ar}

%\today

\maketitle

Let $X$ be a normal scheme and $F$ a coherent sheaf on $X$.
The {\it reflexive hull} or {\it double dual} of $F$ is the sheaf
$F^{**}:=\Hom_X(\Hom_X(F,\o_X),\o_X)$.
The natural map $F\to F^{**}$ kills the torsion
subsheaf of $F$ and the support of $\coker[F\to F^{**}]$
has codimension $\geq 2$. This establishes a functor
from the category of coherent sheaves on $X$ to the
category of reflexive coherent sheaves on $X$.

The same construction works as long as $X$ satisfies
Serre's condition $S_2$, but
otherwise  the double dual is not  $S_2$. One can, however, define
a natural functor from the category of quasi coherent sheaves on $X$ to the
category of quasi coherent sheaves on $X$ that are $S_2$
as sheaves over their support (\ref{hull.defn}).
We denote it by $F\to F^{[**]}$ and call it the {\it  hull} 
 of $F$.

The main question we address in this section is the behaviour
of the hull in families. The motivating example
is the following theorem 
which describes all possible base changes
that create a flat sheaf out of a non-flat sheaf.

\begin{thm}[Flattening decomposition theorem] \cite[Lecture 8]{mumf66}
\label{flat.dec.mumf.thm}
Let $f:X\to S$ be a projective morphism and $F$ a coherent sheaf
on $X$. Then
$S$ can be written as a disjoint union of locally closed subschemes
$S_i\to S$ such that for any $g:T\to S$
the following are equivalent:
\begin{enumerate}
\item  the pull back of
$F$ to $X\times_ST$ is flat over $T$, and
\item  $g$ factors
through the disjoint union $\amalg S_i\to S$.
\end{enumerate}
\end{thm}

Here we study a similar question where instead of
the flatness of $F$ we aim to understand the
flatness of the hulls of the fiber-wise
restrictions $\bigl(F|_{X_s}\bigr)^{[**]}$. 
Note that even in very nice situations, for instance when
$f:X\to S$ is smooth and the restrictions $F|_{X_s}$ are all
torsion free, the  hulls $\bigl(F|_{X_s}\bigr)^{**}$
do not form a sheaf on $X$. Thus the flattening decomposition theorem
does not apply directly.

The  main result  (\ref{hull.exists.thm}) is a close analog of
the flattening decomposition theorem for projective morphisms.
he next formulation is somewhat vague;
see (\ref{hull.exists.thm}) for the precise version.

\begin{thm}\label{intro.thm2}
Let $f:X\to S$ be a projective morphism and $F$ a coherent sheaf
on $X$.
%Assume that the set of points where $F$ is flat over $S$
% is dense in every fiber.
 Then
$S$ can be written as a disjoint union of locally closed subschemes
$S_i\to S$ such that for any $g:T\to S$ the following are equivalent:
\begin{enumerate}
\item  The  hulls
$\{F_t^{[**]}: t\in T\}$ form a flat sheaf on 
$X\times_ST\to T$.
\item  The map $g$ factors
through the disjoint union $\amalg S_i\to S$.
\end{enumerate}
\end{thm}

The cases when  the hull of each $F_t^{\otimes m}$
 is locally free of rank 1 for some $m>0$
has been treated in \cite{hacking, abr-hass}.
In moduli theory, the main application of (\ref{intro.thm2})
is to the hulls $\omega_X^{[m]}$
of $\omega_X^{\otimes m}$, see (\ref{omega.m.comm.cor}, \ref{omega.all.comm.cor}).
  As a consequence we obtain a
well defined  theory of those deformations where
the hulls $\omega_X^{[m]}$ form a flat family.

The first  step in the proof is the construction of
the moduli space of {\it husks} and {\it quotient husks}. 
Generalizing the notion of a hull,
a husk of a quasi coherent sheaf $F$ is a map $q:F\to G$
where $G$ is torsion free over its support
and $q$ is an isomorphism at the generic points.
There are no instability problems for coherent husks, and we prove
that they have a fine moduli space (\ref{q-husk.exists.thm}).
The necessary techniques are taken from \cite{FGA-4, le-p1, le-p2}
with very little change.
Similar ideas have been used in
\cite{honsen, ale-knu, pan-tho, rydh}.

 It is then not difficult 
to identify the hulls among the husks to obtain (\ref{hull.exists.thm}).

The appendix discusses how to extend these results from
projective morphisms to proper morphisms of algebraic spaces.
This is joint work with M.\ Lieblich.

It would also be of interest to obtain a local version of the 
flattening decomposition theorem,
 but our methods are very much global in nature.

\subsection*{Husks}
\medskip

\begin{defn} 
 We say that a quasi coherent sheaf $F$ on a scheme $X$ is
{\it pure} or {\it torsion free} over its support if every associated
prime of $F$ has dimension $=\dim\supp F$, that is, the maximum of the
dimensions of the supports of local sections of $F$.  
In particular,  $\supp F$ is pure dimensional.
We also say that $F$ is pure of dimension $n:=\dim\supp F$.
If $j:X\to Y$ is finite and $F$ is pure
then $j_*F$ is also pure.

For a quasi coherent sheaf $G$ on a scheme $X$,
let $\tors G\subset G$\index{??} denote the subsheaf of those local sections
whose support has dimension $<\dim\supp G$. Thus
$G/\tors G$ is pure.

 Let $f:X\to S$ be a morphism and $F$ a quasi coherent 
sheaf on $X$. We say that $F$ is {\it pure} over $S$  or that it  has  
{\it torsion free fibers} over their support
if for every $s\in S$, the restriction
$F_s$ is pure of the same dimension.

\end{defn}

\begin{defn}\label{husk.field.defn}
  Let $X$ be a scheme over a field $k$, $F$ a quasi coherent 
sheaf on $X$ and  $n:=\dim \supp F$. 
A {\it husk} of $F$ is a  quasi coherent  sheaf $G$ together with a
homomorphism $q:F\to G$ such that
\begin{enumerate}
\item $G$ is pure  of dimension $n$, and
\item $q:F\to G$ is an isomorphism at all $n$-dimensional points of $X$.
\end{enumerate}

If $h\in \ann(F)$ then $h\cdot F=0$, hence $h\cdot G\subset G$
is supported in dimension $<n$, hence $0$. Thus $G$ is also
an $\o_X/\ann(F)$ sheaf and so the particular choice of $X$
matters very little.

Assume that $X$ is projective and $H$ is ample on $X$.  
As in (\ref{hilb.say}), for a coherent sheaf $M$ on $X$ write
$$
\chi \bigl(X, M(tH)\bigr)=: \sum a_i(M)t^i.
\eqno{(\ref{husk.field.defn}.3)}
$$
Set $n:=\dim\supp F$ and let $G$ be a husk of $F$.  
If $F,G$ are coherent then, by (\ref{hilb.say}.1)
$$
a_n(G)=a_n\bigl(F/\tors F\bigr) \qtq{and}
a_{n-1}(G)\geq a_{n-1}\bigl(F/\tors F\bigr).
\eqno{(\ref{husk.field.defn}.4)}
$$
\end{defn}

\begin{say}[Universal husk]\label{univ.husk.say}

The smallest husk of $F$ is $F/\tors F$.
 
There is also a largest or  {\it universal husk} $U(F)$
which can be constructed as follows. Let $R$ be the total ring of
quotients of $\o_X/\ann(F)$. That is, we invert every element
that is a unit at every $n$-dimensional generic point of
 $\o_X/\ann(F)$. Then $F\to F\otimes_XR$ is the universal husk.

First, $F\otimes_XR$ is an $R$-sheaf, hence its associated 
primes are the $n$-dimensional generic points of
 $\o_X/\ann(F)$. By construction, $F\to F\otimes_XR$
is an isomorphism at every
$n$-dimensional generic point of $\supp F$.

Second, let $F\to G$ be any other husk.
Then we get $F\otimes_XR\to G\otimes_XR$ which is an
isomorphism at  every
$n$-dimensional generic point of $\supp F$ hence an isomorphism
of $R$-sheaves. Thus $F\to F\otimes_XR$
factors as $F\to G\to F\otimes_XR$. Since $G$ has no
lower dimensional associated primes,
$G\to F\otimes_XR$ is an injection. Hence the husks of $F$
are the quasi coherent subsheaves of $U(F)$ that contain
$F/\tors F$.

If $n\geq 1$ then  $U(F)$ is never coherent, but it is
the union of coherent husks.
Thus a coherent sheaf has many different coherent husks
and there is  no universal coherent husk.
\end{say}

\begin{defn}\label{husk.defn}
  Let $f:X\to S$ be a morphism and $F$ a quasi coherent 
sheaf on $X$. %which is generically flat on every fiber of $\supp F\to S$.
 Let $n$ be the relative dimension of
$\supp F\to S$.  
A {\it husk} of $F$ is a quasi coherent  sheaf $G$ together with a
homomorphism $q:F\to G$ such that
\begin{enumerate}
\item $G$ is flat and pure over $S$, 
%\item $G_s:=G|_{X_s}$ has no associated primes of dimension $<n$
% for every $s\in S$, and
\item $q:F\to G$ is an isomorphism at every
$n$-dimensional  point of $X_s\cap\supp F$  for every $s\in S$.
Equivalently, if $q_s:F_s\to G_s$ is a husk for every $s\in S$.
\end{enumerate} 
Note that the notion of a husk does depend on $f$.

As before,  $G$ is also
an $\o_X/\ann(F)$ sheaf and so  $X$
matters very little.

Husks are preserved by base change.
That is, if $g:T\to S$ is a morphism, $X_T:=X\times_ST$ and
$g_X:X_T\to X$ the first projection then
$g_X^*q:g_X^*F\to g_X^*G$ is also a husk.
\end{defn}

\begin{lem}\label{husk.funct.lem}
  Let $f:X\to S$ be a morphism and $F$ a quasi coherent 
sheaf on $X$. %which is generically flat on every fiber of $\supp F\to S$.
Let $q:F\to G$ be a husk of $F$.
\begin{enumerate}
\item Let $g:X\to Z$ be a finite $S$-morphism.
Then $g_*G$ is a husk of $g_*F$.
\item Let $h:Y\to X$ be a flat morphism.
Then $h^*G$ is a husk of $h^*F$.
\end{enumerate}
\end{lem}

Proof.  If $g$ is a finite morphism and $M$ is a sheaf
then the associated primes
of $g_*M$ are the images of the associated primes of $M$.
This implies (1). Similarly, if $h$ is flat then
 the associated primes
of $h^*M$ are the preimages of the associated primes of $M$,
 implying (2).
\qed

\begin{defn}\label{HUSK.FUNCT.EFB}  Let $f:X\to S$ be a morphism, 
and  $F$  a  coherent sheaf on $X$.
%which is
%generically flat on every fiber of $\supp F\to S$.
Let ${\it Husk}(F)(*)$ be the functor that to a
scheme $g:T\to S$ associates the
set of all  coherent husks of $g_X^*F$ with proper support over $T$,
where $g_X:T\times_SX\to X$ is the projection.

Let $f:X\to S$ be a projective morphism, 
$H$ an $f$-ample divisor and  $p(t)$  a polynomial. 
Let ${\it Husk}_p(F)(*)$ be the functor that to a
scheme $g:T\to S$ associates the
set of all  coherent husks of $g_X^*F$
with Hilbert polynomial $p(t)$.
\end{defn}

\begin{defn}
 Let $f:X\to S$ be a morphism and $F$ a quasi coherent sheaf
on $X$. The husk of a quotient of $F$ is called a {\it quotient husk} of $F$.
Equivalently, a  quotient husk of $F$ is a  quasi coherent  sheaf 
$G$ together with a
homomorphism $q:F\to G$ such that
\begin{enumerate}%\setcounter{enumi}{4}
\item $G$ is pure over $S$, say of relative dimension $m$ and
\item $q:F\to G$ is surjective at all $m$-dimensional points of $X_s\cap\supp G$
for every $s\in S$. 
\end{enumerate} 

As in (\ref{HUSK.FUNCT.EFB}), %for a  coherent sheaf $F$ 
${\it QHusk}_p(F)(*)$ denotes the functor that to a
scheme $g:T\to S$ associates the
set of all  coherent quotient husks of $g_X^*F$
with Hilbert polynomial $p(t)$,
where $g_X:T\times_SX\to X$ is the projection.

\end{defn}

The first  existence theorem  is the following.

\begin{thm}\label{q-husk.exists.thm}
  Let $f:X\to S$ be a projective morphism, 
$H$ an $f$-ample divisor,  $p(t)$  a polynomial and 
 $F$  a  coherent sheaf on $X$.
Then  ${\it QHusk}_p(F)$ is bounded, proper,  separated
and it has a fine moduli space ${\rm QHusk}_p(F)$.
\end{thm}

(The construction establishes ${\rm QHusk}_p(F)$
as an algebraic space. There does not seem to be any obvious
ample line bundle on it, so its projectivity may be
a subtle question.)

Proof.  We  start by establishing the valuative criterion
of properness and separatedness.
%%REFERENCE  TO OTHER SECTION \ref??)
 Then we check
that ${\it QHusk}_p(F)$ is bounded. 
The moduli space ${\rm QHusk}_p(F)$ is then constructed
using   the theory of quotients by algebraic group actions.

As a preliminary step, note that the problem is local, thus
we may assume that $S$ is affine. Then $f,X,F$ are defined over
a finitely generated subalgebra of $\o_S$, hence we may assume
in the sequel that $S$ is of finite type.

\medskip

\ref{q-husk.exists.thm}.1 
{\it The valuative criterion of  separatedness and properness.}

Let $T$ be the spectrum of a DVR with closed point $0\in T$
and generic point $t\in T$. Given $g:T\to S$
we have $g_X^*F$ 
%which, by assumption, is flat
%outside $Z_T:=g_X^{-1}(Z)$ 
where $g_X:T\times_SX\to X$ is the projection
%and $Z$ 
as in (\ref{genflat.defn}).
 By assumption, we also have a
husk $q_t:F_t\to G_t$; set $F'_t:=\im q_t$. By (\ref{Quot-schemes.say}),
 there is a unique flat
quotient $g_X^*F\onto F'$ that agrees with $F'_t$ over the generic point.

Further, there is a closed subset $B_t\subset \supp G_t$
such that $\dim B_t< \dim\supp G_t$ and 
$F'_t\to G_t$  is an isomorphism
outside $B_t$.
Let $B_T\subset X_T$ be the closure of $B_t$.

$F'$ 
is flat over $T$,
its generic fiber is pure
and its special fiber is  pure
outside a subset $Z_0\subset X_0$ of dimension $<\dim\supp G_0$.
Furthermore, 
$G_t$ and $F'$  are naturally isomorphic over $X_t\setminus  B_t$.
Thus we can glue them to get a single sheaf
$G'$ defined on  $X_T\setminus (Z_0\cup B_0)$.
By construction,  $G'$ is flat and pure over $T$.

Let $j:X_T\setminus (Z_0\cup B_0)\into X_T$ be the injection.
Set $G:=j_* G'$.
Since $Z_0\cup B_0$ has codimension $\geq 2$
in $\supp G_T$,
 the push forward $G$ is coherent by (\ref{push.finite.lem}).
The fibers of $G$   are pure
 by (\ref{pff.char.lem}.4). In particular, $G$ is
flat over $T$ and  $G_0$ is a husk of $F_0$.
By (\ref{pff.char.lem}.4),  $G$ is the only extension of $G'$ 
that is pure over $T$. These show that
${\it QHusk}(F)$ satisfies the 
 valuative criterion of  separatedness and properness.

Furthermore,  $G_0$ has
 the same Hilbert polynomial as $G_t$.
\medskip

\ref{q-husk.exists.thm}.2 {\it Boundedness.}

 $F$ is $m(F)$-regular (\ref{m-reg.say}) for some  $m(F)$.
 We show that 
all  quotient husks
$q:F\to G$ are $m$-regular  for some  $m$
depending only on $m(F)$ and $p(t)$.
This is   a problem  on the individual fibers,
so from now on we assume that $X\subset \p^n$ is a projective
scheme over a field. We also use that this assertion holds for
all quotients of $F$ (\ref{Quot-schemes.say}).

The proof is by induction on $\deg p(t)=\dim\supp G$.

If $\dim \supp G=0$ then any $m$ works.

Assume next that $\dim \supp G=1$. 
Let $G'\subset G$ be the image of $F$.
The Hilbert polynomial of $G'$ is $p(t)-c$ where $c$ is the length of $G/G'$. 
Pick any $K''\subset K:=\ker q$ such that
$K/K''$ has length $c$.
Then $G'':=F/K''$ is a quotient of $F$ with  Hilbert polynomial
$p(t)$, hence $m$-regular for some $m$
depending only on $m(F)$ and $p(t)$.
Since $G''\to G$ has 0-dimensional kernel and cokernel,
 this implies that $G$ is also $m$-regular.

Assume now that $\dim\supp G=n\geq 2$. 
After a field extension, we may assume that the base field is infinite.
Let $H$ be a general hyperplane.
Then  $F|_H$ is also $m$-regular and 
$F|_H\to G|_H$ is a  
quotient husk with Hilbert polynomial $p(t)-p(t-1)$ by (\ref{husk.funct.lem.3}).
Hence, by induction, $G|_H$ is $m_1$-regular for some  $m_1$
depending on $m(F)$ and $p(t)$.
$G$ is then $m$-regular for some  $m$
depending on $m(F)$ and $p(t)$ by (\ref{inductive.mreg.lem}).
\medskip

\ref{q-husk.exists.thm}.3 {\it Construction of ${\rm QHusk}_p(F)$.}

The existence of ${\rm QHusk}_p(F)$ is a local problem on
$S$.
As we noted in (\ref{husk.defn}), we can replace
$X$ with any larger scheme.
Thus we may assume that $X=\p^n_S$.

By boundedness,  we can choose $m$ such that
for any quotient husk $F\to G$ with Hilbert polynomial $p(t)$, 
$G_s(m)$ is generated by global sections and
its higher cohomologies vanish. Thus each
$G_s(m)$ can be written as a quotient of
$\o_{X_s}^{\oplus p(m)}$.

As in (\ref{Quot-schemes.say}), let 
$$
Q_{p(t)}:=\quot_{p(t)}^0(\o_X^{\oplus p(m)})
\subset \quot(\o_X^{\oplus p(m)})
$$
be the universal family of quotients  $q_s:\o_{X_s}^{\oplus p(m)}\onto M_s$
that have Hilbert polynomial $p(t)$,  are pure,
have no higher cohomologies and
the induced map
$$
q_s:H^0\bigl(X_s ,  \o_{X_s}^{\oplus p(m)}\bigr)
\to H^0\bigl(X_s ,M_s\bigr)
$$
is an isomorphism.

Let $\pi:Q_{p(t)}\to S$ be the structure map, 
 $\pi_X:Q_{p(t)}\times_S X\to X$ the second projection and 
$M$  the universal sheaf on $Q_{p(t)}\times_SX$.

By (\ref{hom.sch.exists}) there is an  open subscheme 
 $W_{p(t)}=\uhom^0 (\pi_X^*F,M)\subset \uhom (\pi_X^*F,M)$
parametrizing  those maps from $\pi_X^*F$ to $M$
that are surjective outside a subset of dimension $\leq n-1$.
Let  $\sigma:W_{p(t)}\to Q_{p(t)}$ be the structure map, and 
$\sigma_X:W_{p(t)}\times_S X\to Q_{p(t)}\times_S X$  the
 fiber product.

Note that $W_{p(t)}$ parametrizes triples
$$
w:=\Bigl[F_w\stackrel{r_w}{\to} G_w\stackrel{q_w}{\twoheadleftarrow}
\o_{X_w}(-m)^{\oplus p(m)}\bigr)
\Bigr]
$$
where $r_w:F_w\to G_w$ is a quotient husk with Hilbert polynomial $p(t)$
and $q_w(m):\o_{X_w}^{\oplus p(m)}\to G_w(m)$
is a surjection that induces an isomorphism on the spaces
of global sections.

Let $w'\in W_{p(t)}$ be another point corresponding to the triple
$$
w':=\Bigl[F_{w'}\stackrel{r_{w'}}{\to} G_{w'}
\stackrel{q_{w'}}{\twoheadleftarrow}
\o_{X_{w'}}(-m)^{\oplus p(m)}\bigr)
\Bigr].
$$
such that 
$$
\bigl[F_w\stackrel{r_w}{\to} G_w\bigr]\cong 
\bigl[F_{w'}\stackrel{r_{w'}}{\to} G_{w'}\bigr].
$$
Then the difference between $w$ and $w'$
comes from the different ways that
we can write $G_w\cong G_{w'}$  as quotients of 
$\o_{X_w}(-m)^{\oplus p(t)}$.
Since we assume that  
$q_w(m), q_{w'}(m)$ induce isomorphisms on the spaces of global sections,
% H^0\bigl(X_w ,  \o_{X_w}^{\oplus p(m)}\bigr)
%\rightrightarrows H^0\bigl(X_w ,G_w(m)\bigr)=H^0\bigl(X_w ,G_{w'}(m)\bigr)
%$$
%are  isomorphisms, 
the different choices of
$q_w$ and $q_{w'}$  correspond to  different bases
in $H^0\bigl(X_w ,G_w(m)\bigr)$. Thus the
fiber of 
$$
\mor(*,W_{p(t)})\to {\it QHusk}_p(F)(*)
\qtq{over} 
\pi\circ\sigma (w)=\pi\circ\sigma (w')=:s\in S
$$ is a
principal homogeneous space under
the group scheme
$$
GL(p(t), k(s))=\aut\Bigl(H^0\bigl(X_s ,G_s(m)\bigr)\Bigr).
$$
 Let $G$ be the group scheme
$GL(p(t),S)$. Then 
$G$ acts on $W_{p(t)}$ and
${\rm QHusk}_p(F)=W_{p(t)}/G$  \cite{k-quot, ke-mo}. 
\qed
%%REF TO MODULI SECTION   by (\ref??).  

\begin{defn}\label{genflat.defn} 
 Let $f:X\to S$ be a morphism and $F$ a quasi coherent sheaf
on $X$.  Let $n=\max_{s\in S}\dim \supp \bigl(F|_{X_s}\bigr)$
We say that $F$ is {\it generically flat} on every fiber of $\supp F\to S$ if
$F$ is flat at every $n$-dimensional point of every fiber $X_s$.
If $F$ is coherent, then this is equivalent to the following:

There is a subscheme $Z\subset X$ such that
\begin{enumerate}
\item $F|_{X\setminus Z}$ is flat over $S$, and
\item $\dim (X_s\cap Z)<n$ for every $s\in S$.
\end{enumerate}
%If $f$ is proper and $F$ is coherent, then $\supp F\to S$ is also proper
%and, by flatness,  $s\mapsto \dim (X_s\cap \supp F)$
%is locally constant. To simplify notation we always assume
%that it is actually constant.
\end{defn}

\begin{cor}\label{husk.exists.thm}
Notation and assumptions as in (\ref{q-husk.exists.thm}).
\begin{enumerate}
\item  ${\it Husk}_p(F)$ is bounded,  separated
and it has a fine moduli space ${\rm Husk}_p(F)$
which is an open subspace of ${\rm QHusk}_p(F)$.
\item Assume that  $F$  is generically flat
 on every fiber of $\supp F\to S$.
Then ${\it Husk}_p(F)$ is proper and 
${\rm Husk}_p(F)\subset {\rm QHusk}_p(F)$ is  closed.
\end{enumerate}
\end{cor}

Proof. Let $q_U:\pi^*F\to G_U$ be the universal
quotient husk over $\pi:{\rm QHusk}_p(F)\to S$.
For a surjective map with flat target, it
is an open condition to
be fiber-wise isomorphic (cf.\ \cite[22.5]{mats-cr}). 
Thus there is a closed subset 
$Z\subset {\rm QHusk}_p(F)\times_SX$
such that $q_U:\pi^*F\to G_U$ is a fiberwise isomorphism
exactly outside $Z$. Then ${\it Husk}_p(F)\subset {\rm QHusk}_p(F)$
is the largest open subset  over which the fiber dimension
of $Z\to {\rm QHusk}_p(F)$ is less than $\deg p(t)$. 

Assume next that  $F$  is generically flat
 on every fiber of $\supp F\to S$. Then there is a closed
subscheme $W\subset {\rm QHusk}_p(F)\times_SX$
such that the  fiber dimension
of $W\to {\rm QHusk}_p(F)$ is less than $\deg p(t)$, and 
$F$ is flat and $q$ is surjective outside $W$.
Then $\ker q$ is also flat outside $W$, hence
$\dim \supp\ker q_s<n$ is  a closed condition.
Note that $q_s:F_s\to G_s$ is a husk iff 
$q_s$ is generically injective, that is iff $\dim \supp \ker q_s<n$.
This proves (2).\qed

\begin{say}[Restriction of husks]\label{husk.funct.lem.3}
Let $q:F\to G$ be a husk or a quotient husk 
and $H\subset X$   a Cartier divisor.
When is $F|_H\to G|_H$ a husk or a quotient husk?

First of all, we need that  $G|_H$ be flat and pure of dimension  $(n-1)$.
The first of these holds if $H$ does not contain any
associated prime of any $G_s$,   cf.\ \cite[Thm.22.5]{mats-cr}.
Since $\bigl(G|_H\bigr)|_s=G_s|_H$, the second condition is satisfied
if $H$ does not contain any
associated prime of any of the hulls $G_s^{[**]}$; see
(\ref{S_2.no.ext.lem}) and   
(\ref{hull.field.defn}).

Second, for any $s\in S$, we have an exact sequence
$$
0\to K_s\to F_s \to G_s\to G_s/F_s\to 0.
$$
Tensoring with $\o_H$ is exact if
$\tor^1(F_s/K_s, \o_H)=0$ and
$\tor^1(G_s/F_s, \o_H)=0$.
Every associated prime of $F_s/K_s$
is an associated prime of $G_s$.
Thus both vanishings hold if 
$H$ does not contain any associated prime of $G_s$ or of $G_s/F_s$
for every $s\in S$.

In particular, if $G$ is coherent and the residue fields of
$S$ are infinite, then these conditions hold for
general members of any base point free linear system
on $X$. (See (\ref{hulls.of.coh}) for the required coherence of $G_s^{[**]}$.)
\end{say}

\subsection*{Hulls}
\medskip

\begin{defn}\label{hull.field.defn}
  Let $X$ be a scheme over a field $k$ and $F$ a quasi coherent 
sheaf on $X$.  Set $n:=\dim\supp F$.  
A husk  $q:F\to G$ is called  {\it tight} if
$q$ is onto  at all $(n-1)$-dimensional
 points of $X$. 
There is a  unique maximal tight husk
$q:F\to F^{[**]}$, called the {\it hull} 
(or $S_2$-hull, see (\ref{hull.over.vars})) of $F$.

$F^{[**]}$, as a subset of the universal husk $U(F)$ 
defined in (\ref{univ.husk.say}),
 is generated by all local sections $\phi\in U(F)$ such that
$\phi$ is a local section of $q(F)$ at all $(n-1)$-dimensional
 points of $\supp G$. 

If $X$ itself is normal, $F$ is coherent   and $\supp F=X$, then
$F^{[**]}$ is the usual reflexive hull $F^{**}$ of $F$.

The hull of a nonzero sheaf is also nonzero, in contrast with the
reflexive hull which kills all torsion sheaves.
\end{defn}

\begin{lem}\label{hull.over.vars}
 Let $X$ be a scheme over a field $k$ and $F$ a  coherent 
sheaf on $X$. 
\begin{enumerate} 
\item Let $r:F\to G$ be any tight husk.  Then $q:F\to F^{[**]}$
extends uniquely to an injection $q_G:G\into F^{[**]}$.
\item $F^{[**]}$ is the unique tight husk which is $S_2$ over its support.
\item $F^{[**]}$ is the smallest $S_2$ husk of $F$. That is, if
 $r:F\to G$ is any husk such that $G$ is $S_2$ over its support, then
$r$ factors as  $F\to F^{[**]}\into G$.
\item  Let $Z\subset X$ be  a closed subset    such that
$\dim Z\leq n-2$  and 
$F/\tors F$ is $S_2$  over $X\setminus Z$.
Let $j:X\setminus Z\to X$ denote the injection.
Then 
$$
F^{[**]}=j_*\bigl( (F/\tors F)|_{X\setminus Z}\bigr).
$$
\item Assume that $X$ is projective, $H$ is ample on $X$, $F\to G$ is any 
coherent husk
and $n=\dim\supp F$. Then, 
\begin{enumerate} 
\item $a_n\bigl(F^{[**]}\bigr)=a_n\bigl(F/\tors F\bigr)$ and
$a_{n-1}\bigl(F^{[**]}\bigr)=a_{n-1}\bigl(F/\tors F\bigr)$.
%\eqno{(\ref{hulls.of.coh}.1)}
\item $a_n\bigl(F^{[**]}\bigr)=a_n\bigl(G\bigr)$ and
$a_{n-1}\bigl(F^{[**]}\bigr)\leq a_{n-1}\bigl(G\bigr)$,
%\eqno{(\ref{hulls.of.coh}.2)}
\item  equality holds iff $G\subset F^{[**]}$.
\end{enumerate}
Thus the hull minimizes $a_{n-1}$ and maximizes the rest
of the Hilbert polynomial.
\end{enumerate}
\end{lem}

Proof. The first property holds by definition.

Let $r:F\to G$ be a tight husk such that
$G\subsetneq F^{[**]}$. Pick any $\phi\in F^{[**]}\setminus G$ and
a function $f\in \o_X$ which is invertible at
all $n$-dimensional generic points of $\supp F$ such that $f\phi\in G$.
Then $f\phi\in G/fG$ has $\leq (n-2)$-dimensional support,
thus $G$ is not $S_2$.

Conversely, with $f$ as above, 
let $\phi\in F^{[**]}/fF^{[**]}$ be a local section which has
$\leq (n-2)$-dimensional support.
Then $\langle F^{[**]},f^{-1}\phi\rangle \in U(F)$
is also a tight husk of $F$. Thus $\phi\in   fF^{[**]}$
and so $F^{[**]}/fF^{[**]}$ has no nonzero local sections with
$\leq (n-2)$-dimensional support. Thus $F^{[**]}$ is $S_2$, hence (2) holds.

Let $r:F\to G\subset U(F)$ be a husk which is $S_2$.
Pick any local section $\phi\in F^{[**]}$. 
Then $\langle G, \phi\rangle/G$ is supported in
dimension $\leq (n-2)$. Since $G$ is $S_2$, this implies that
$\phi\in G$, proving (3).

 (4) is discussed in greater detail in (\ref{pff.char.lem}).

 (5.a) follows from (\ref{hilb.say}.1) and, together with
(\ref{husk.field.defn}.4), it implies (5.b).
If $a_n\bigl(F^{[**]}\bigr)=a_n\bigl(G\bigr)$ and
$a_{n-1}\bigl(F^{[**]}\bigr)= a_{n-1}\bigl(G\bigr)$, then, by
(\ref{hilb.say}.2), $F\to G$ is a tight husk, hence
$G\subset F^{[**]}$ by (1).
 \qed

\begin{say}[Hulls of coherent sheaves]\label{hulls.of.coh}
The hull $F^{[**]}$ of a coherent sheaf $F$
is almost always coherent.
For instance, this holds if $X$ is of finite type over a field
or over an excellent ring.

To see this, we can assume that $X$ is affine
and replace $F$ by $F/\tors F$. Then there is
 a sequence of
subsheaves $0=F_0\subset \dots\subset F_n=F$ such that
every $F_{m+1}/F_m$ is isomorphic to an ideal sheaf in  $\o_{X_m}$ for some
integral subscheme $X_m\subset X$ of dimension $n$.

By (\ref{hull.over.vars}.4),  $F\to F^{[**]}$ is left exact
on sequences of pure $n$-dimensional sheaves. Thus
it is sufficient to prove
that  the hull of any  ideal sheaf is coherent.
In turn this follows if the hull of  $\o_{X_m}$ is coherent.
By (\ref{hull.over.vars}.3) the hull of  $\o_{X_m}$ is contained in the
 normalization of  $\o_{X_m}$. Thus
 $\o_{X_m}^{[**]}$ is coherent whenever  the
normalization is coherent.
\end{say}

\begin{defn}\label{hull.defn}
  Let $f:X\to S$ be a morphism and $F$ a quasi coherent 
sheaf.
 Let $n$ be the relative dimension of
$\supp F\to S$.  
A {\it hull} of $F$ is a  husk $q:F\to G$
such that, for every $s\in S$,
 the induced map 
$q_s:F_s\to G_s$ is a hull (\ref{hull.field.defn}).

We see in (\ref{hulls.unique}) that a hull is unique if it exists.
Note that  if a hull exists then
$F$ is generically flat on every fiber of $\supp F\to S$.

It is clear from the definition that hulls are preserved by base change.
That is, if $g:T\to S$ is a morphism, $X_T:=X\times_ST$ and
$g_X:X_T\to X$ the first projection then
$g_X^*q:g_X^*F\to g_X^*G$ is also a hull.
\end{defn}

\begin{lem}\label{hulls.unique} 
 Let $f:X\to S$ be a morphism of finite type and $F$ a  coherent 
sheaf on $X$.
%which is generically flat on every fiber of $\supp F\to S$.
 Let $n$ be the relative dimension of
$\supp F\to S$. 
\begin{enumerate}
\item Let  $q:F\to G$ be a hull
and set $Z:=\supp G/F$. Then  $G$ is coherent, 
$\dim (X_s\cap Z)\leq n-2$ for every $s\in S$ and 
$F/\tors F$ is flat over $X\setminus Z$.

\item 
Conversely, let $Z\subset X$ be any closed subset such that
$\dim (X_s\cap Z)\leq n-2$ for every $s\in S$ and 
$F/\tors F$ is flat over $X\setminus Z$.
Let $j:X\setminus Z\to X$ denote the injection.
Then 
$$
G=j_*\bigl( (F/\tors F)|_{X\setminus Z}\bigr).
$$
In particular, $F$ has at most one hull.
\end{enumerate}
\end{lem}

Proof. $G_s=\bigl(F|_s\bigr)^{[**]}$ is coherent by (\ref{hulls.of.coh}),
thus $G$ is coherent by the Nakayama lemma.
The rest of the first part is clear from the definition.
To see the converse,
let  $q:F\to  G$ be a hull and
$Z$ any closed subset such that
$\dim (X_s\cap Z)\leq n-2$ for every $s\in S$ and $Z\supset\supp \coker q$. 
Then $F/\tors F$ and $G$ are
isomorphic over $X\setminus Z$, hence
$F/\tors F$ is flat over $X\setminus Z$.
Furthermore, by (\ref{pff.char.lem}),  
$$
G=j_*\bigl( G|_{X\setminus Z}\bigr)=
 j_*\bigl( (F/\tors F)|_{X\setminus Z}\bigr).
$$
Thus $G$ is unique.\qed

\begin{defn}  Let $f:X\to S$ be a projective morphism and 
 $F$  a  coherent sheaf on $X$. 
%which is
%generically flat on every fiber of $\supp F\to S$.
For a scheme $g:T\to S$ set
 ${\it Hull}(F)(T)=1$ if $g_X^*F$ has a hull
and ${\it Hull}(F)(T)=\emptyset$ if $g_X^*F$ does not have a hull,
where $g_X:T\times_SX\to X$ is the projection.
\end{defn}

\begin{defn} A morphism  $g:\bar S\to S$
is a {\it locally closed decomposition} of $S$
if 
\begin{enumerate}
\item for every connected component $\bar S_i\subset \bar S$,
the restriction of $g$ to  $\bar S_i$ is a locally closed embedding, and
\item $g$ is one-to-one and onto on geometric points.
\end{enumerate}
\end{defn}

The second existence  theorem  is the following.

\begin{thm}[Flattening decomposition for hulls]\label{hull.exists.thm}
  Let $f:X\to S$ be a projective morphism and 
 $F$  a  coherent sheaf on $X$.
%which is generically flat
% on every fiber of $\supp F\to S$.
Then 
\begin{enumerate}
\item ${\it Hull}(F)$ is bounded,  separated
and it has a fine moduli space ${\rm Hull}(F)$.
\item  The  structure map 
 ${\rm Hull}(F)\to S$ is a locally closed decomposition.
\end{enumerate}
\end{thm}

Proof. We construct the  locally closed decomposition
${\rm Hull}(F)\to S$ by first identifying a closed stratum
and then using  induction.

Let $n$ be the maximal fiber dimension of $\supp F\to S$.

For any  point $s\in S$ write
$$
\chi \bigl(X_s, (F_s)^{[**]}(tH)\bigr)=:p_s(t)=:
a_n(s)t^n+a_{n-1}(s)t^{n-1}+O(t^{n-2}).
$$
By (\ref{flat.dec.mumf.thm}) and 
(\ref{up-lo.s.cont.prop}.3), $\{p_s(t):s\in S\}$ is a finite set of polynomials.
Let $p(t)=a_nt^n+a_{n-1}t^{n-1}+O(t^{n-2})$ be the polynomial which
lexicographically maximizes the triple
$\bigl(a_n(s), -a_{n-1}(s), p_s\bigr)$ for all $s\in S$.
(Note the minus sign before $a_{n-1}(s)$.)
\medskip

(\ref{hull.exists.thm}.3) {\it Claim.} Every quotient husk of $F$
with Hilbert polynomial $p(t)$ is a hull.
\medskip

Proof. This holds after any base change, but, for simplicity of notation,
we work directly over $S$.

Let $F\to G$ be a  quotient husk
with Hilbert polynomial $p(t)$. 
The following exact sequences define $K$ and $F'$:
$$
0\to F'\to G\to G/F'\to 0
\qtq{and}
0\to K\to F\to F'\to 0.
$$
Since $G$ is flat over $S$ and
the fiber dimension of $\supp(G/F')\to S$ is less than $n$,
see see that $F'$ is flat over $S$ at the generic points of its support
in each fiber.

Therefore, for every $s\in S$,
$K_s\to F_s\to F'_s\to 0$ is exact at all generic points of $\supp G_s$.
Thus $a_n(F_s)\geq a_n(F'_s)$. On the other hand, we assumed that
$a_n(F'_s)=a_n(G_s)$ is the largest possible. Thus
$a_n(F_s)= a_n(F'_s)$
and so  $\supp K\to S$ has fiber dimension $<n$ over $S$.
In particular,    $F\to G$ is  a  husk.

Since $F_s\to G_s$ is a husk, $a_{n-1}(F_s)\leq a_{n-1}(G_s)$ by
(\ref{hull.over.vars}.5).
 By our choice,  $a_{n-1}(G_s)$ is the smallest possible,
hence  $a_{n-1}(F_s)= a_{n-1}(G_s)$ and so
$G_s\subset F_s^{[**]}$ by (\ref{hull.over.vars}.5).
Since $p(t)$ is maximized, this implies that $G_s= F_s^{[**]}$.
Thus $F\to G$ is a hull. \qed

By (\ref{q-husk.exists.thm}),  ${\rm QHusk}_p(F)\to S$ is proper.
As we proved, it parametrizes hulls, hence
${\rm QHusk}_p(F)\to S$ is a monomorphism 
(\ref{monom.defn}, \ref{hulls.unique}.2).
A proper monomorphism is a closed embedding
(\ref{monom.defn}); let $S_p\subset S$ denote the image
of ${\rm QHusk}_p(F)\to S$.

We can now replace $S$ by $S\setminus S_p$ and
conclude by induction on the cardinality of
$\{p_s(t):s\in S\}$.\qed

\begin{defn}[Monomorphisms]\label{monom.defn} A morphism of schemes
$f:X\to Y$ is a {\it monomorphism} if for every scheme $Z$
the induced map of sets
$\mor(Z,X)\to \mor(Z,Y)$ is an injection.

By \cite[IV.17.2.6]{EGA} this is equivalent to assuming that
$f$ be universally injective and unramified.

A closed or open embedding is a monomorphism.
Other typical example of monomorphisms is the 
 normalization of the node with a point missing, that is
$\a^1\setminus\{-1\}\to (y^2=x^3+x^2)$ given by
$(t\mapsto (t^2-1, t^3-t)$.

A proper monomorphism $f:Y\to X$ is a closed embedding.
Indeed, a proper monomorphism is injective
on geometric points, hence finite. 
Thus it is a closed embedding
iff $\o_X\to f_*\o_Y$ is onto. By the Nakayama lemma this is
equivalent to
$f_x:f^{-1}(x)\to x$ being an isomorphism for every
 $x\in f(Y)$. By passing to geometric points,
we are down to the case when 
 $X=\spec k$, $k$ is algebraically closed and
  $Y=\spec A$ where $A$ is an Artin $k$-algebra. 

If $A\neq k$ then there are at least 2 different $k$ maps
$A\to k[\epsilon]$, thus $\spec A\to \spec k$ is
not a monomorphism.
\end{defn}

D.\ Rydh pointed out that, besides hulls, it is also of interest to 
consider $F/\tors F$, which is the smallest husk of $F$.
In this case, the method of (\ref{hull.exists.thm}) gives the following
if we first maximize $a_n$ and then minimize $p(t)$.

\begin{prop}[Flattening decomposition for pure quotients]
Let $f:X\to S$ be a projective morphism and $F$ a coherent sheaf
on $X$.
%  which is generically flat
% on every fiber of $\supp F\to S$. 
Then $S$ can be written as a disjoint union of locally closed subschemes
$S_i\to S$ such that for any $g:T\to S$
the following are equivalent:
\begin{enumerate}
\item   $g_X^*F/\tors g_X^*F$  is flat and pure.
\item  $g$ factors
through the disjoint union $\amalg S_i\to S$.
\end{enumerate}
\end{prop}

\subsection*{Applications}

{\ }\medskip

Applying (\ref{hull.exists.thm}) to the relative dualizing sheaf
gives the following result.

\begin{cor} \label{omega.m.comm.cor}
Let $f:X\to S$ be  projective and equidimensional.
Let $Z\subset X$  be  a closed subscheme such that
 $\codim (X_s, Z\cap X_s)\geq 2$
 for every $s\in S$ and  $(X\setminus Z)\to S$ is flat with Gorenstein fibers.
Then, for any $m$ there is a locally closed decomposition $S_m\to S$
such that for any  $g:T\to S$ the following are equivalent
\begin{enumerate}
\item  
$\omega_{X\times_ST/T}^{[m]}$
is flat over $T$ and commutes with base change.
\item  $g$ factors through $S_m\to S$.
\end{enumerate}
\end{cor}

Proof. The question is local on $S$, thus we may assume that
there is a finite surjection $\pi:X\to \p^n_S$. 
One can now define  $\omega_{X/S}$ as
$$
\omega_{X/S}:=\Hom_{\p^n_S}\bigl(\pi_*\o_X, \omega_{\p^n_S/S}\bigr).
$$
In general,  $\omega_{X/S}$ does not commute with base change
but, by assumption,  its restriction to $X\setminus Z$ is locally free.

We claim that $S_m={\rm Hull}\bigl(\omega_{X/S}^{\otimes m}\bigr)$.

Given $g:T\to S$, let $j_T:X\times_ST\setminus Z\times_ST\to X\times_ST$
be the inclusion. Then 
$$
\omega_{X\times_ST/T}^{[m]}=\bigl(j_T\bigr)_*
g_X^*\omega_{X\setminus Z/S}^{\otimes m}.
$$
If $T\mapsto \omega_{X\times_ST/T}^{[m]}$
 commutes with restrictions to
the fibers of $X\times_ST\to T$, then 
$\omega_{X\times_ST/T}^{[m]}$ has $S_2$ fibers, hence
$\omega_{X\times_ST/T}^{[m]}$
 is the hull of $\omega_{X\times_ST/T}^{\otimes m}$.

Conversely, by (\ref{hulls.unique}), 
if $\omega_{X\times_ST/T}^{\otimes m}$ has a hull
then it is $\omega_{X\times_ST/T}^{[m]}$ and
it commutes with further base changes by (\ref{hull.defn}).\qed

\begin{cor} \label{omega.all.comm.cor}
Let $f:X\to S$ be  projective and equidimensional.
Let $Z\subset X$  be  a closed subscheme such that
 $\codim (X_s, Z\cap X_s)\geq 2$
 for every $s\in S$ and  $(X\setminus Z)\to S$ is flat with Gorenstein fibers.
Assume in addition that there is an $N>0$ such that
$\omega_{X_s}^{[N]}$ is locally free for every $s\in S$.

Then there is a locally closed decomposition $S^*\to S$
such that a morphism $g:T\to S$ factors through $S^*$
iff
$\omega_{X\times_ST/T}^{[m]}$
is flat over $T$ and commutes with base change
for every $m\in \z$.
\end{cor}

Proof. 
Let $S_i\to S$ be as in (\ref{omega.m.comm.cor}).
Take $S^*$ to be the fiber product of  the morphisms 
$S_1\to S, \dots, S_N\to S$.\qed
\medskip

\subsection*{Hilbert polynomials of non-flat sheaves}
\medskip

\begin{say}[Hilbert polynomials]
\label{hilb.say}

Let $X$ be a projective scheme of dimension $n$ and $H$ an ample divisor.
For a coherent sheaf $F$, write its  Hilbert polynomial as
$$
\chi(X,F(tH))=a_n(F)t^n+\cdots+a_0(F).
$$
Then $a_n(F)\geq 0$ and $a_n(F)= 0$ iff $\dim\supp F<n$.

(\ref{hilb.say}.1)  Let $u:F\to G$ be a map of coherent sheaves which is an
isomorphism outside a subset $Z\subset X$ of dimension
$\leq n-r$. Then 
 $a_i(F)=a_i(G)$ for $n\geq i> n-r$ and
if $\dim\supp\ker u<n-r$ then $a_{n-r}(G)\geq a_{n-r}(F)$.
Indeed note that
$$
\chi(X,F(tH))-\chi(X,G(tH))=\chi(X,(\ker u)(tH))-\chi(X,(\coker u)(tH)).
$$
By assumption,  both $\coker u$ and $\ker u$ are
supported on $Z$, hence their Hilbert polynomials have
degree $\leq n-r$.

(\ref{hilb.say}.2) Conversely, let  $u:F\to G$ be a map of sheaves which is an
isomorphism at the generic points.
If  $a_i(F)=a_i(G)$ for $n\geq i> n-r$ and
$\dim\supp\tors F\leq n-r$, then
$u$ is an
isomorphism outside a subset of dimension
$\leq n-r$.

(\ref{hilb.say}.3) Let $f:X\to S$ be a projective morphism of 
pure relative dimension $n$.
Let $F$ be a sheaf on $X$.
Fix an integer $r$ and assume that there is  a closed subscheme 
$Z\subset X$ such that
$F$ is flat over $X\setminus Z$ and
 $\dim_s (Z\cap X_s)\leq n-r$  for every $s\in S$.

Pick  an $f$-very ample divisor  $H$.
Every $s\in S$ has an open neighborhood
$U$ such that for general  $H_0,\dots, H_{n-r}\in |H_{f^{-1}(U)}|$,
the restriction
$F|_{H_0\cap\cdots\cap H_{n-r}}$ is flat over $U$.
In particular, the Hilbert polynomial
$\chi\bigl(H_0\cap\cdots\cap H_{n-r}, F(m)\bigr)$ is well defined.

For each $u\in U$, the  Hilbert polynomial
$\chi\bigl(H_0\cap\cdots\cap H_{n-r}\cap X_u, F(m)\bigr)$
determines the top $r$ coefficients of the
Hilbert polynomials  $\chi(X_u, F(t)|_{X_u})$. 
Thus we conclude the following.

(\ref{hilb.say}.4)
 Under the 
 assumptions of (\ref{hilb.say}.3), the top $r$ coefficients of the
Hilbert polynomials  $\chi(X_s, F(t)|_{X_s})$ are
locally constant on $S$.
\end{say}

\begin{defn}
 Let $p_1$ and $p_2$  be two polynomials.
We say that $p_1\leq p_2$ if $p_1(t)\leq p_2(t)$ for all $t\gg 0$.

For example, if $F_1\subset F_2$ are coherent sheaves on 
a projective scheme  and $p_i$ is the Hilbert
polynomial of
$F_i$ then $p_1\leq p_2$ and equality holds iff $F_1=F_2$.
\end{defn}

\begin{prop}\label{up-lo.s.cont.prop}  Let $f:X\to S$ be a projective morphism 
 and $H$ an $f$-ample Cartier divisor.
Let $F$ be a coherent sheaf on $X$.
For a point $s\in S$, set $F_s:=F|_{X_s}$.
\begin{enumerate}
\item  The  Hilbert polynomial function
$$
s\mapsto  \chi \bigl(X_s, F_s(tH)\bigr)
$$
is constructible and upper semi continuous on $S$. 
\item  Assume that  $F$ is  generically flat on every fiber of
 $\supp F\to S$.
Then the  Hilbert polynomial
$$
s\mapsto  \chi \bigl(X_s, (F_s/\tors F_s)(tH)\bigr)
$$
is constructible and lower semi continuous on $S$. 
\item  If  $F$ is  generically flat on every fiber of
 $\supp F\to S$ and  $\dim\supp \tors F_s\leq \dim\supp  F_s-2$ 
for every $s\in S$
 then the  Hilbert polynomial of the hull
$$
s\mapsto  \chi \bigl(X_s, (F_s)^{[**]}(tH)\bigr)
$$
is  constructible and upper semi continuous on $S$. 
\end{enumerate}
\end{prop}

Proof. We may assume that $S$ is reduced.
 By generic flatness \cite[Lec.8]{mumf66}
 there is a dense open subset $S^0\subset S$ such that $F$ 
and $F/\tors F$ are
flat over $S^0$ and $F/\tors F$ has pure fibers.
Thus  $\chi \bigl(X_s, F_s(tH)\bigr)$ and 
$\chi \bigl(X_s, (F_s/\tors F_s)(tH)\bigr)$
are both locally constant on $S^0$. By Noetherian induction
we conclude that both functions are  constructible.

It is thus enough to check  semi continuity
when $S$ is the spectrum of a DVR.
Let $0\in S$ be the closed point and $g\in S$ the general point.
Let $\tors_0F\subset F$ be the torsion supported on $X_0$.
Then $F/\tors_0F$ is flat over $S$ and so
$$
\chi \bigl(X_g, F_g(tH)\bigr)=
\chi \bigl(X_0, (F/\tors_0F)\otimes \o_{X_0}(tH)\bigr).
$$
There is an exact sequence
$$
(\tors_0F)\otimes \o_{X_0}\to F\otimes \o_{X_0}\to
(F/\tors_0F)\otimes \o_{X_0}\to 0,
$$
hence
$$
\chi \bigl(X_0, F\otimes \o_{X_0}(tH)\bigr)
\geq \chi \bigl(X_0, (F/\tors_0F)\otimes \o_{X_0}(tH)\bigr)=
\chi \bigl(X_g, F_g(tH)\bigr),
$$
which proves upper semi continuity in the first case.

If  $F$ is  generically flat on every fiber of
 $\supp F\to S$, then 
$$
\dim \supp (\tors_0F)\otimes \o_{X_0}<\dim \supp F_0.
$$
Thus $(\tors_0F)\otimes \o_{X_0}$ maps to the torsion
subsheaf of $F_0$ and
$F_0/\tors F_0$ is also a quotient of 
$(F/\tors_0F)\otimes \o_{X_0}$. Therefore
$$
\chi \bigl(X_0, (F_0/\tors F_0)(tH)\bigr)
\leq \chi \bigl(X_0, (F/\tors_0F)\otimes \o_{X_0}(tH)\bigr)=
\chi \bigl(X_g, F_g(tH)\bigr),
$$
which proves lower semi continuity in the second case.

In order to prove (\ref{up-lo.s.cont.prop}.3),
let $j:U\into X$ be an open set such that
$F/\tors F$ is torsion free and flat over $U$
and $(X\setminus U)\cap X_{s_g}$ has codimension $\geq 2$
for a generic point $s_g\in S$. There is an  open
neighborhood $U_1$ of $s_g$ such that 
$(X\setminus U)\cap X_s$ has codimension $\geq 2$
for every point $s\in U_1$. Furthermore, 
 $j_*\bigl(F|_U\bigr)$ is flat 
and has $S_2$ fibers over a nonempty open $U_2\subset U_1$.
Thus 
$$
\o_{X_s}\otimes j_*\bigl(F|_U\bigr)\cong \bigl(F|_{X_s}\bigr)^{[**]},
$$
and so $\chi \bigl(X_s, (F_s)^{[**]}(tH)\bigr)$
is locally  constant on $U_2$. As before, this proves constructibility.

As before, it is  enough to check upper semi continuity
when $S$ is the spectrum of a DVR and $F$ is torsion free.
In particular, $F$ is flat over $S$.

Let $j:U\into X$ be an open set such that
$F$ is  flat over $U$,
 $(X\setminus U)\cap X_g$ has codimension $\geq 2$
and $(X\setminus U)\cap X_0$ has codimension $\geq 1$.
Then  $G:=j_*\bigl(F|_U\bigr)$ is flat over $S$ and
$F\to G$
is generically an isomorphism on every fiber.
On the generic fiber, $G_g\cong (F_g)^{[**]}$.
On the special fiber we know that $G_0$ is torsion free and there is a map
$F_0\to G_0$ which is an isomorphism at all generic points.

Since $F$ and $G$ are both flat over $S$,
$$
\chi \bigl(X_0, F_0(tH)\bigr)=\chi \bigl(X_g, F_g(tH)\bigr)\qtq{and}
\chi \bigl(X_0, G_0(tH)\bigr)=\chi \bigl(X_g, G_g(tH)\bigr).
$$
Furthermore, since $F_g$ is torsion free,
$$
\deg_t  \Bigl(\chi \bigl(X_g, G_g(tH)\bigr)-\chi \bigl(X_g, F_g(tH)\bigr)\Bigr)
\leq n-2,
$$
hence also
$$
\deg_t  \Bigl(\chi \bigl(X_0, G_0(tH)\bigr)-\chi \bigl(X_0, F_0(tH)\bigr)\Bigr)
\leq n-2.
$$
By assumption, $\dim\supp \tors F_0\leq n-2$,
hence $F_0\to G_0$ is an  isomorphism at all
codimension 1 points by (\ref{hilb.say}.2). Hence by 
(\ref{hull.over.vars}.1)
there is an injection $G_0\into (F_0)^{[**]}$.
Thus
$$
\chi \bigl(X_0, (F_0)^{[**]}(tH)\bigr)
\geq \chi \bigl(X_0, G_0(tH)\bigr)=
\chi \bigl(X_g, G_g(tH)\bigr)=\chi \bigl(X_g, (F_g)^{[**]}(tH)\bigr.\qed
$$

\medskip

\begin{exmp} The condition on $\dim\supp\tors F_s$ in
(\ref{up-lo.s.cont.prop}.3)
 is necessary. Let $C,S$ be smooth projective curves
and $X=C\times S$. Let $f:X\to S$ be the projection and
 $F$  the ideal sheaf of a point $(c,s_0)\in X$. 
For $s\neq s_0$, $F_s\cong \o_C$ has Hilbert polynomial
$t\deg H+1-g(C)$.  On the other hand,
$F_{s_0}\cong \o_C(-c)+k(c)$ and its hull is
$\o_C(-c)$ with Hilbert polynomial
$t\deg H-g(C)$.
\end{exmp}

\subsection*{Quot-schemes}

{\ }\medskip

\begin{say}[Quot-schemes]\label{Quot-schemes.say} \cite{FGA-4}
 Let $f:X\to S$ be a morphism and $F$ a  coherent sheaf
on $X$.  ${\it Quot}(F)(*)$ denotes the functor that to a
scheme $g:T\to S$ associates the
set of all   quotients of $g_X^*F$
that are flat over $T$ with proper support,
where $g_X:T\times_SX\to X$ is the projection.

If $F=\o_X$, then a quotient can be identified with a subscheme
of $X$, thus ${\it Quot}(\o_X)={\it Hilb(X)}$, the Hilbert functor.

If $H$ is an $f$-ample divisor and  $p(t)$  a polynomial,
then ${\it Quot}_p(F)(*)$ denotes those flat quotients that have
Hilbert polynomial $p(t)$.

By \cite{FGA-4},  ${\it Quot}_p(F)$ is bounded, proper,  separated
and it has a fine moduli space ${\rm Quot}_p(F)$.
See \cite[Sec.4.4]{sernesi} for a detailed proof.
If $F=\o_X$, then ${\rm Quot}(\o_X)=\hilb(X)$, the Hilbert scheme of $X$. 

Note that one can write $F$ as a quotient of $\o_{\p^n}(-m)^{r}$
for some $m,r$, thus ${\it Quot}_p(F)$ can be viewed as a
subscheme of ${\it Quot}(\o_{\p^n}^{r})$.
The theory of  ${\it Quot}(\o_{\p^n}^{r})$ is essentialy the same
as the study of the Hilbert functor, discussed in
\cite{mumf66} and \cite[Sec.I.1]{rc-book}.

\end{say}

\begin{say}[Castelnuovo-Mumford regularity] \label{m-reg.say}
Let $F$ be a coherent
sheaf on $\p^n$. We say that $F$ is {\it $m$-regular}
if $H^i(\p^n,F(m-i))=0$ for $i\geq 1$.
See   \cite[Sec.I.1.8]{laz-book}
for a detailed treatment.

It is known that if $F$ is $m$-regular then it is also $m'$-regular for
every $m'\geq m$ and the multiplication maps
$$
H^0(\p^n, F(m'))\otimes H^0(\p^n, \o_{\p^n}(1))  
\to  H^0(\p^n, F(m'+1))
$$
are surjective for $m'\geq m$.
Thus $F(m)$ is generated by global sections and so $F$ is a quotient
of the sum of $h^0(\p^n, F(m))$ copies of $\o_{\p^n}(-m)$.
In particular, all $m$-regular sheaves with Hilbert polynomial
$p(t)$ are quotients of $\o_{\p^n}(-m)^{p(m)}$, hence
they form a  bounded family.

One can almost get a uniform vanishing theorem
for $H^0(\p^n,F(-r))$ as follows.
Let $F$ be a coherent  sheaf on $\p^n$
and $s\in H^0(\p^n, F)$ a section whose support has dimension $d$. Then
$$
h^0(\p^n, F(m))\geq h^0(\p^d, \o_{\p^d}(m))=\binom{m+d}{d}.
$$
In particular,  if $F$ is $m$-regular and
$h^0(\p^n,F(m))=\chi(\p^n,F(m))=:r$ then every
section of $F(m-r)$ has 0-dimensional support.

\end{say}

\begin{lem}\label{inductive.mreg.lem} Let $G$ be a coherent sheaf
on $\p^n$  with Hilbert polynomial $p(t)$.
Assume that $G$ has no associated primes of dimension $<2$.
Let $H\subset \p$ be a hyperplane that does not contain any of the
 associated primes of $G$ and assume that $G|_H$ is $m_1$-regular.

Then $G$ is $m$ regular for some $m$ depending only on
$m_1$ and $p(t)$.
\end{lem}

Proof. Using the cohomology sequence of
$$
0\to G(r-1)\to G(r)\to G|_H(r)\to 0
$$
we conclude that $H^i(X, G(r-1))\cong H^i(X, G(r))$ for
$i\geq 2$ and $r\geq m_1-i+1$. Thus, by Serre's vanishing,
$H^i(X, G(r))=0$ for $ i\geq 2$ and $r\geq m_1-i$. 

For $i=1$ we have only an exact sequence
$$
H^0(X, G(r))\stackrel{b(r)}{\to} 
H^0(X\cap H, G|_H(r))\to H^1(X, G(r-1))\stackrel{c(r)}{\to} 
 H^1(X, G(r))\to 0, 
$$
which shows that $b(r)$ is onto iff $c(r)$ is an isomorphism.

We also have a commutative square
$$
\begin{array}{ccc}
H^0(X, G(r))\otimes H^0(\p^n, \o_X(1)) & \to & H^0(X, G(r+1))\\*[4pt]
b(r)\otimes \sigma\downarrow\hphantom{b(r)\otimes q} && 
\hphantom{b(r+1)}\downarrow b(r+1)\\
H^0(X\cap H, G|_H(r))\otimes H^0(H, \o_H(1)) & 
\stackrel{t(r)}{\to}  & H^0(X\cap H, G|_H(r+1))
\end{array}
$$
where $\sigma:H^0(\p^n, \o_X(1))\to H^0(H, \o_H(1))$
is the (surjective)  restriction.
Since $G|_H$ is $m_1$-regular, $t(r)$ is onto for $r\geq m_1$ (\ref{m-reg.say}).

This shows that if $b(r)$ is onto for some $r\geq m_1$
then $b(r+1)$ is also onto.
Thus, if $b(r)$ is onto 
then  $b(s)$ is onto for every $s\geq r$ and
$c(s)$ is  an isomorphism for every $s\geq r$.
Again  by Serre's vanishing, this gives that
$H^1(X, G(r))=0$.

Otherwise
$H^1(X, G(r))\neq 0$ but then
$h^1(X, G(r))> h^1(X, G(r+1))$. In either case
we get that
$$
H^1(X, G(r))=0\qtq{for $r\geq m_1+h^1(X, G(m_1))$.}
$$
Since $h^1(X, G(m_1))=h^0(X, G(m_1))-p(m_1)$, we are done
if we can bound $h^0(X, G(m_1))$ from above.
Since $G$ has no 0-dimensional associated primes,
$$
h^0(X, G(m_1))\leq \sum_{i\geq 0} h^0(H, G|_H(m_1-i)),
$$
and the latter sum is finite by (\ref{m-reg.say}) 
and bounded by induction if
$G|_H$ has no 0-dimensional associated primes.
The latter follows from our assumptions.
\qed
\medskip

The following is proved in 
\cite[III.7.7.8--9]{EGA}, see also \cite[4.6.2.1]{la-mb} and \cite[2.1.3]{lieb}.

\begin{defn-lem}\label{hom.sch.exists}
 Let $f:X\to S$ be proper. Let $F,L$ be coherent sheaves
on $X$ such  that $L$ is flat over $S$. Then there is an $S$-scheme
$\uhom(F,L)$ parametrizing homomorphisms from $F$ to $L$.
That is, for any $g:T\to S$, there is a natural isomorphism
$$
\Hom_T(g_X^*F, g_X^*L)\cong \mor_S\bigl(T,  \uhom(F,L)\bigr),
$$ 
where $g_X: T\times_SX\to X$ is the fiber product
of $g$ with the identity of $X$.
\end{defn-lem}

Proof. Note that there is a natural identification
between
\begin{enumerate}
\item homomorphisms $\phi:F\to L$, and
\item quotients  $\Phi:(F+L)\onto M$ which  induce an isomorphism 
$\Phi|_L:L\cong M$.
\end{enumerate}

Let  $\pi:\quot(F+L)\to S$ denote the quot-scheme parametrizing 
quotients of $F+L$
with universal quotient
$u:\pi_X^*(F+L)\to M$, where
$\pi_X$  denotes the induced map  
$\pi_X:\quot(F+L)\times_SX\to X$.

Consider now the restriction of $u$ to
$u_L:\pi_X^*L\to M$. By the Nakayama lemma, for a map between sheaves
 it is an open condition to
be surjective. For a surjective map with flat target, it
is an open condition to
be fiber-wise injective (cf.\ \cite[22.5]{mats-cr}). 
Thus there is an open subset
$$
 \quot^0(F+L)\subset \quot(F+L)
$$
which parametrizes those quotients
$v:F+L\to M$ which induce an isomorphism $v_L:L\cong M$.
Thus $\uhom(F,L)=\quot^0(F+L)$.\qed

\subsection*{Push forward and $S_2$}

{\ }\medskip

Here we collect some well known results
about normalization, pushing forward and $S_2$-sheaves.

\begin{lem}\label{S_2.no.ext.lem}
 Let $R$ be a Noetherian ring and $M$ an $R$-module.
Assume that each associated prime of $M$ has dimension $n$.
The following are equivalent:
\begin{enumerate}
\item If $r\in R$ is not contained in any associated prime of $M$
then every associated prime of $M/rM$ has dimension $(n-1)$.
\item If $N\supset M$ has the same associated primes as $M$
and every associated prime of $N/M$ has dimension $\leq (n-2)$
then $N=M$.
\end{enumerate}
\end{lem}

Proof. To see (2) $\Rightarrow$ (1), 
pick $n\in N$ such that the associated primes of $Rn/Rn\cap M$
have dimension $\leq (n-2)$. There is an $r\in R$ which
is not contained in any associated prime of $Rn\cap M$ such that $rn\in M$.
Then $rn$ leads to an associated prime of $M/rM$ of dimension
$\leq (n-2)$. By (1) we obtain that $rn\in rM$ and so
$n\in M$. 

Conversely, assume that there is a submodule $rM\subset M'\subset M$
such that every associated prime of $M'/rM$ has dimension
$\leq (n-2)$. Then $N:=r^{-1}M'\supset M$ shows that $r^{-1}M'=M$,
Thus $M'=rM$, which gives  (1) $\Rightarrow$ (2).\qed

\medskip

The proof of the next lemma is essentially the same
as the coherence argument in (\ref{hulls.of.coh}).

\begin{lem}\label{push.finite.lem}
  \cite[IV.5.11.1]{EGA}  Let $X$ be the spectrum of a 
Nagata ring
$R$,
$i:U\to X$
 the immersion of an open set and 
 $F$  a coherent sheaf on $U$.
Then $i_*F$ is coherent iff  $\codim(   \bar x\cap (X\setminus U), \bar x)\geq
2$
for every  associated prime $x\in X$ of $F$.
\end{lem}

\begin{lem} \label{pff.char.lem} Let $X$ be an affine scheme, $i:U\to X$ 
an open  immersion  and $W=X\setminus U$. 
Let $F$ be a coherent sheaf on $U$ and $G$ a quasi coherent sheaf on $X$
 such that $G|_U\cong F$. The
following  conditions  are equivalent:
\begin{enumerate}
\item  $G\cong i_*F$.
\item  For $a,b\in I_W$, if  the sequence
$$
0\to F\stackrel{(b,-a)}{\longrightarrow} F+F
\stackrel{(a,b)}{\longrightarrow} F
$$
is exact then so is
$$
0\to G\stackrel{(b,-a)}{\longrightarrow} G+G
\stackrel{(a,b)}{\longrightarrow} G.
$$
\item $\depth_W(G)\geq 2$.
\item If $a\in I_W$ does not vanish on any associated
prime of $F$ then  $G/aG$ has no associated
prime supported on $W$.
\item  $H^i_W(X,G)=0$ for $i=0,1$.
\item  $H^0_W(X,G)=0$ and for every coherent sheaf $Q$ such that $\supp
Q\subset W$ every extension
$$
0\to G\to G'\to Q\to 0\qtq{splits. }
$$
\end{enumerate}
\end{lem}

Proof.  We first prove that (\ref{pff.char.lem}.1) and 
(\ref{pff.char.lem}.2) are equivalent. Since $F$ is
coherent, it has only finitely many associated primes. Choose $a\in I_W$ which
is not contained in any associated prime of $F$ and 
$b\in I_W$ which is not contained in any associated prime of $F/aF$. 

Then
$0\to F\stackrel{(b,-a)}{\longrightarrow} F+F
\stackrel{(a,b)}{\longrightarrow} F$ is exact.
Since $i_*$ is left exact, this implies (\ref{pff.char.lem}.2). 
Conversely, assume (\ref{pff.char.lem}.2).
  $G|_U\cong F$, thus there is a natural homomorphism $G\to i_*F$. Let $K\subset
G$ be its kernel. Every element of $K$ is killed by  a power of $I_W$, thus if
$K\neq 0$ then we can choose $0\neq g\in K$ such that $I_W\cdot g=0$,
 in particular
$ag=bg=0$. This is impossible since (\ref{pff.char.lem}.2) is left exact.

Thus $G\to i_*F$ is an injection; let $C$ be its cokernel. 
Every element of $C$ is killed by  a power of $I_W$, thus if
$C\neq 0$ then we can choose $g'\in i_*F\setminus G$ such that 
$ag',bg'\in G$. $(bg',-ag')\in G+G$ is in the kernel of the map $(a,b)$, 
hence by
exactness there is a $g\in G$ such that $bg=bg'$ and $ag=ag'$. 
 This implies that
$g-g'=0$, a contradiction.

(\ref{pff.char.lem}.2) implies that $(a,b)$ is a $G$-sequence of length two, 
hence
$\depth_W(G)\geq 2$. Conversely, if $\depth_W(G)\geq 1$ then none of the
associated primes of $G$ are contained in $W$.
Every other associated prime of $G$ is also an associated prime
of $F$, hence $G$ has only finitely many associated primes.
Therefore  one can  choose $a\in I_W$
which is not contained in any associated prime of $F$. 
If $\depth_W(G)\geq 2$ then none of the
associated primes of $G/aG$ are contained in $W$, hence one can  choose $b\in
I_W$ which is not contained in any associated prime of $G/aG$. This shows that
(\ref{pff.char.lem}.2) is equivalent to (\ref{pff.char.lem}.3). 

(\ref{pff.char.lem}.4) is a restatement of
(\ref{pff.char.lem}.3).

For any quasi coherent sheaf $G$  there is an exact sequence
$$
0\to H^0_W(X,G)\to H^0(X,G)\to H^0(U, G|_U)\to H^1_W(X,G)\to H^1(X,G).
$$
Thus $H^0_W(X,G)=0$ iff $G\to i_*(G|_U)$ is an injection.
Since $X$ is affine, $H^1(X,G)=0$ thus
$H^1_W(X,G)=0$ iff $G\to i_*(G|_U)$ is a surjection. These show that
 (\ref{pff.char.lem}.1) and
(\ref{pff.char.lem}.5) are equivalent.

If $H^0_W(X,G)=0$ then $G\to i_*F$ is an injection. Thus if $G\neq i_*F$ then
we have a nonsplit extension. If $G\to G'$ is as in (\ref{pff.char.lem}.6)
 then $G'|U\cong F$
 gives a homomorphism $G'\to i_*F$ which is a  splitting of $G\to G'$
if $G\cong i_*F$.
\qed
\medskip

%\newpage
\subsection*{Appendix: Algebraic space case}

(joint with M.\ Lieblich)

We consider the case when $S$ is a Noetherian algebraic space
and $f:X\to S$ is a proper morphism of algebraic spaces.

\begin{say}[Flat families of coherent sheaves]\label{flat.stack.say}
Let $f:X\to S$ be a proper morphism. The functor of flat families of
coherent sheaves ${\it Flat}(X/S)$ is represented by
an algebraic stack  ${\rm Flat}(X/S)$ which is locally of finite type
but very nonseparated; cf.\ \cite[4.6.2.1]{la-mb}.

However,  as in (\ref{q-husk.exists.thm}.1), 
${\rm Flat}(X/S)$ satisfies the existence part of the 
valuative criterion
of properness. That is, if $T$ is the spectrum of a DVR with generic point
$t_g$ then every morphism $t_g\to {\rm Flat}(X/S)$
extends to $T\to {\rm Flat}(X/S)$.

In fact, an even stronger property holds:
\medskip

\ref {flat.stack.say}.1 {\it Theorem.} \cite[5.2.2, 5.7.9]{ray-gru}
Let $U/S$ be a separated $S$-space of finite type and
$g:U\to  {\rm Flat}(X/S)$ a morphism. Then there is an
$S$-space $\bar U\supset U$ which is proper over $S$
such that $g$ extends to a morphism $\bar g:\bar U \to  {\rm Flat}(X/S)$.
\medskip

\end{say}

\begin{say}[Construction of ${\rm QHusk}(F)$]

Let $\sigma:{\rm Flat}(X/S)\to S$ be the structure morphism and
let $U_{X/S}$ denote the universal family over
$ {\rm Flat}(X/S)$. There is an open substack
$$
{\rm Flat}^{n}(X/S)\subset {\rm Flat}(X/S)
$$
parametrizing those sheaves that are pure of dimension $n$. 
Let $U^n_{X/S}$ be the corresponding   universal family.

Consider $X\times_S {\rm Flat}^n(X/S)$ with coordinate projections
$\pi_1,\pi_2$. The stack
$$
\uhom\bigl(\pi_1^*F, \pi_2^* U^n_{X/S}\bigr)
$$
parametrizes all maps from the sheaves $F_s$ to sheaves $N_s$
 that are pure of dimension $n$ (\ref{hom.sch.exists}). 

We claim that ${\rm QHusk}(F)$ is an open substack
of $\uhom\bigl(\pi_1^*F, \pi_2^* U^n_{X/S}\bigr)$.
Indeed, as in the proof of
(\ref{hom.sch.exists}), 
for a map of sheaves $M\to N$ with $N$ flat over $S$, it is
an open condition to be surjective  at the generic points
of the support.

As in (\ref{q-husk.exists.thm}.1), 
we see that  ${\rm QHusk}(F)$ is separated.
\end{say}

Putting these together, 
and using that an algebraic  stack 
whose diagonal is a monomorphism is an algebraic space
(see, for instance, \cite[Sec.8]{la-mb}),
we obtain the first
 existence theorem:

\begin{thm}\label{asp.husk.exists.thm}
  Let $f:X\to S$ be a proper morphism of algebraic spaces and 
 $F$  a  coherent sheaf on $X$.
%which is generically flat
% on every fiber of $\supp F\to S$.
Then 
\begin{enumerate}
\item ${\it QHusk}(F)$ is  separated
and it has a fine moduli space ${\rm QHusk}(F)$.
\item Every irreducible component of ${\rm QHusk}(F)$
is proper over $S$. \qed
\end{enumerate}
\end{thm}

\begin{say}[Construction of ${\rm Hull}(F)$]\label{contr.hull.say}

In a flat family of coherent shaves, it is an open condition
to be $S_2$ over their support and for a map to a flat sheaf
it is also
an open condition to be an isomorphism at the codimension 1 points
of their support.
This implies that 
${\rm Hull}(F)$ is an open subspace of
${\rm QHusk}(F)$.  

We claim that ${\rm Hull}(F)$
is of finite type. First, it is locally of finite type
since  ${\rm QHusk}(F)$ is.
Second, we  claim that $\red {\rm Hull}(F)$ is dominated
by an algebraic space of finite type.
In order to see this, consider  the (reduced) structure map 
$\red{\rm Hull}(F)\to \red S$.
It is an isomorphism at the generic points, hence there is an
open dense $S^0\subset \red S$ such that
$S^0$ is isomorphic to an open subspace
of $\red {\rm Hull}(F)$. 
Repeating this for $\red S\setminus S^0$,
by Noetherian induction 
we eventually write $\red {\rm Hull}(F)$ as a disjoint union of
finitely many locally closed subspaces of $\red S$.
(We do not claim, however, that every 
irreducible component of  $\red {\rm Hull}(F)$ is a 
locally closed subspace of $\red S$.)

These together imply that  ${\rm Hull}(F)$
is of finite type. (Indeed, if $U\to V$ is a surjection,
$U$ is of finite type and $V$ is locally of finite type
then $V$ is of finite type.)

As in (\ref{hulls.unique}.2), the structure map
${\rm Hull}(F)\to  S$ is a monomorphism.
However, in the nonprojective case, it need not be
a locally closed decomposition (though we do not know any examples).
We can summarize these considerations in the following theorem.
\end{say}

\begin{thm}[Flattening decomposition for hulls]\label{asp.hull.exists.thm}
  Let $f:X\to S$ be a proper morphism of algebraic spaces and
 $F$  a  coherent sheaf on $X$.
%which is generically flat
% on every fiber of $\supp F\to S$.
Then 
\begin{enumerate}
\item ${\it Hull}(F)$ is   separated
and it has a fine moduli space ${\rm Hull}(F)$.
\item  ${\rm Hull}(F)$ is an algebraic space of finite type over $S$.
\item  The  structure map 
 ${\rm Hull}(F)\to S$ is a surjective monomorphism. \qed
\end{enumerate}
\end{thm}

In some cases one can see that ${\rm Hull}(F)\to S$
is a locally closed decomposition using the following
valuative criterion of locally closed embeddings.

\begin{prop}\label{lcemb.val.crit}
 Let $f:X\to Y$ be a morphism of finite type.
Then $f$ is a locally closed embedding iff
\begin{enumerate}
\item $f$ is a monomorphism, and
\item if $T$ is the spectrum of a DVR and $g:T\to Y$ a morphism
such that $g(T)\subset f(X)$ then $g$ lifts to
$g_X:T\to X$.
\end{enumerate}
\end{prop}

Proof. Since $f$  is a monomorphism, it is quasi-finite.
Take any proper $\tilde f:\tilde X\to Y$ extending $f$ and
then its Stein factorization.
We obtain an algebraic space $\bar X\subset X$ and a finite morphism
$\bar f:\bar X\to Y$ extending $f$.
Set $Z:=\bar X\setminus X$.
If $Z=\bar f^{-1}\bar f (Z)$ then
$$
f(X)=\bar f\bigl(\bar X\bigr)\setminus \bar f(Z),
$$
and
$\bar f(Z)\subset \bar f\bigl(\bar X\bigr)\subset Y$ are closed embeddings.
Thus $f(X)\subset Y$ is locally closed and 
$f:X\to f(X)$ is a proper monomorphism hence an isomorphism
by (\ref{monom.defn}).

Otherwise, there are points $z\in Z$  and
$x\in X$ such that $\bar f(z)=\bar f(x)$. Let
$T$ be the spectrum of a DVR and $h:T\to \bar X$ a morphism
which maps the closed point to $z$ and the generic point
to $X$. Set $g:=\bar f\circ h$. Then
$g(T)\subset f(X)$ and the only lifting of
$g$ to $T\to \bar X$ is $h$, but $h(T)\not\subset X$.\qed

\begin{say}[Proof of (\ref{hull.exists.thm})]

One can get another proof of (\ref{hull.exists.thm}) using
(\ref{asp.hull.exists.thm}) and (\ref{lcemb.val.crit})
as follows.

Since $f:X\to S$ is projective, there is an
$f$-ample divisor $H$ and we can decompose
${\rm Hull}(F)=\amalg_p{\rm Hull}_p(F)$
according to the Hilbert polynomials.
In order to prove
that each ${\rm Hull}_p(F)\to S$ is a locally closed embedding,
we check the valuative criterion (\ref{lcemb.val.crit}).

We have $X_T\to T$ and a coherent sheaf $F_T$
such that the hulls of $F_{t_g}$ and of $F_{t_0}$
have the same Hilbert polynomials.

Let  $F_{t_g}\to G_{t_g}$ be the hull over the generic point and
 extend $G_{t_g}$ to a husk  $F_T\to G_T$.

Let $G_{t_0}^{[**]}$ denote the hull of  $G_{t_0}$. Then
the composite
$F_{t_0}\to G_{t_0}\to G_{t_0}^{[**]}$
is the hull of $F_{t_0}$.
By assumption and  by flatness
$$
\chi\bigl(X_0, G_{t_0}^{[**]}(tH)\bigr)=p(t)
=\chi\bigl(X_0, G_{t_g}(tH)\bigr)=
\chi\bigl(X_0, G_{t_0}(tH)\bigr).
$$
Hence, by (\ref{hilb.say}.2), $G_{t_0}= G_{t_0}^{[**]}$
and so  $G_{t_0}= G_{t_0}^{[**]}$ is the hull of $F_{t_0}$. 
Thus $G_T$ defines the lifting $T\to {\rm Hull}_p(F)$. \qed
\end{say}

 \begin{ack} I thank D.\ Abramovich, J.\ Lipman and D.\ Rydh 
for many useful comments and corrections.
Partial financial support  was provided by  the NSF under grant number 
DMS-0500198.
\end{ack}

\bibliography{refs}

\vskip1cm

\noindent Princeton University, Princeton NJ 08544-1000

\begin{verbatim}kollar@math.princeton.edu\end{verbatim}

\end{document}